\documentclass[a4paper,12pt]{article}
\usepackage{amsmath, amsthm, amsfonts, amsxtra, amssymb, psfrag, graphicx, color}

\newtheorem{definition}{Definition}
\newtheorem{theorem}[definition]{Theorem}
\newtheorem*{conjecture}{Conjecture}
\newtheorem{proposition}[definition]{Proposition}
\newtheorem{corollary}[definition]{Corollary}
\newtheorem{lemma}[definition]{Lemma}
\graphicspath{{figures/}}

\def\ni{\noindent}

\def\proof{\ni{\sl Proof.}\ }
\def\qed{\hfill\fbox{\hbox{}}\medskip}

\def\P{\mathcal{P}}

\title{Diametral Pairs of Linear Extensions}
\author{Graham Brightwell\footnote{Department of Mathematics, London School of Economics, London, UK. E-mail: G.R.Brightwell@lse.ac.uk} \ and Mareike Massow\footnote{Institut f\"{u}r Mathematik, Technische Universit\"{a}t Berlin, Berlin, Germany. E-mail: massow@math.tu-berlin.de. Research supported by the Research Training Group ``Methods for Discrete Structures" (DFG-GRK 1408) and the Berlin Mathematical School.}}

\begin{document}

\maketitle

\begin{abstract}
	Given a finite poset~$\P$, we consider pairs  of linear extensions of~$\P$ with maximal distance, where the distance between two linear extensions~$L_1, L_2$ is the number of pairs of elements of~$\P$ appearing in different orders in~$L_1$ and~$L_2$. A diametral pair maximizes the distance among all pairs of linear extensions of~$\P$. Felsner and Reuter defined the linear extension diameter of~$\P$ as the distance between a diametral pair of linear extensions.

We show that computing the linear extension diameter is NP-complete in general, but can be solved in polynomial time for posets of width~3. 

Felsner and Reuter conjectured that, in every diametral pair, at least one of the linear extensions reverses a critical pair. We construct a counterexample to this conjecture. On the other hand, we show that a slightly stronger property holds for many classes of posets: We call a poset \emph{diametrally reversing} if, in every diametral pair, both linear extensions reverse a critical pair. Among other results we show that interval orders and 3-layer posets are diametrally reversing. From the latter it follows that almost all posets are diametrally reversing.

\end{abstract}

\section{Introduction}

\begin{figure}[ht]
\label{LEG_N}
	\begin{center}
	\psfrag{1}{$1$}
	\psfrag{2}{$2$}
	\psfrag{3}{$3$}
	\psfrag{4}{$4$}
	\psfrag{a}{$a$}
	\psfrag{b}{$b$}
     	\psfrag{a}{$1234$}
      	\psfrag{b}{$2134$} 
      	\psfrag{c}{$1243$}
      	\psfrag{d}{$2143$}
      	\psfrag{e}{$2413$}
	\psfrag{P}{$\mathcal{N}$}
	\psfrag{G(P)}{$G(\mathcal{N})$}
    	\includegraphics[scale=0.4]{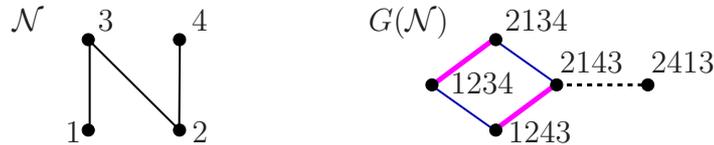}
	\caption{The $\mathcal{N}$ and its LE-graph with the swap colouring.}
  	\end{center}
\end{figure}

Suppose we are given a finite poset~$\P$ and we are interested in the family of its linear extensions (called LEs in the following). If~$\P$ is a chain, then the situation is not very interesting, since there is only one LE, $\P$ itself. If~$\P$ is not a chain, we can ask for a pair of ``maximally different'' LEs of~$\P$. Let us define the \emph{distance} between two LEs~$L_1$ and~$L_2$ as the number of pairs of elements of~$\P$ appearing in different orders in~$L_1$ and~$L_2$, or the number of \emph{reversed pairs} or \emph{reversals} between $L_1$ and $L_2$. We are then interested in \emph{diametral pairs}, which are pairs of LEs of~$\P$ maximising the distance among all pairs of LEs of~$\P$. Such a diametral pair of LEs corresponds to a diametral pair of vertices in the \emph{linear extension graph}~$G(\P)$, where the vertices of~$G(\P)$ are the LEs of~$\P$ and two such vertices are adjacent if the corresponding LEs have distance 1 (see Figure~\ref{LEG_N} for an example). Note that an LE-graph comes with a natural edge colouring, the \emph{swap colouring}, where each edge is coloured with the incomparable pair of elements that is swapped along it. 

LE-graphs were originally defined by Pruesse and Ruskey in \cite{pruesse91}. That article is mainly concerned with the existence of a Hamilton path in the LE-graph; see also \cite{ruskey93} and \cite{west93}. There has been subsequent research on other structural properties of~$G(\P)$; see~\cite{reuter}, \cite{reuter-comp},  \cite{naatz} and \cite{naatz-diss}. Felsner and Reuter~ \cite{felsner-led} define the \emph{linear extension diameter}~$led(\P)$ of a poset~$\P$ as the diameter of~$G(\P)$, or equivalently as the distance between a diametral pair of LEs of~$\P$. In \cite{felsner-led}, some bounds on the LE-diameter are given, and it is determined explicitly for the class of posets called generalized crowns. Also a closed formula for the LE-diameter of the Boolean lattices is conjectured (see Section~\ref{concluding}). 

Is it possible to compute the LE-diameter of a given poset? Given a graph, it is of course easy to determine its diameter. But determining $led(\P)$ in time polynomial in the size of $\P$ may be harder, since the number of LEs of a poset is typically exponentially large. The first main result of this paper is that, given a poset $\P$ and a natural number $k$, it is NP-complete to determine whether $led(\P) \geq k$. However, we show that the LE-diameter of posets of width at most 3 can be computed in polynomial time.

\begin{figure}[ht]
  	\begin{center}
      	\psfrag{x}{$u$}
      	\psfrag{y}{$v$} 
    	\includegraphics[scale=0.35]{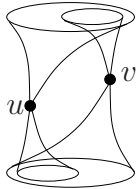}
	\caption{Elements $u$ and $v$ form the critical pair $(u,v)$.}
  	\label{crit_pair}
	\end{center}
\end{figure}

Since we cannot efficiently find a diametral pair of LEs of a given poset, we study properties of  \emph{diametral LEs}, that is, LEs contained in a diametral pair. What makes an LE diametral? A conjecture in~\cite{felsner-led} suggests a connection to \emph{critical pairs}: Two elements $u,v$ of $\P$ form a critical pair $(u,v)$ if they are incomparable in $\P$ and fulfill $\downarrow \!\! u \subseteq \downarrow \!\! v $ and $\uparrow \!\! v \subseteq \uparrow \!\! u$ (see Figure~\ref{crit_pair}). Here, $\downarrow \!\! u$ is short for $\downarrow_\P \!\! u$ and denotes the set of elements smaller than $u$ in $\P$, also called the downset of $u$ or the set of its predecessors. Similarly, $\uparrow_\P \!\! u$ or $\uparrow \!\! u$ denotes the set of elements greater than $u$ in $\P$, which is also called the upset of $u$ or the set of its successors.

Note that critical pairs are ordered, and that they live on the ground set of $\P$ (unlike the diametral pairs). Critical pairs are the relevant incomparable pairs for building a realiser of $\P$ (see e.g.~\cite{trotter92}). It is a classic observation that every poset contains a critical pair -- with the exception of chains, because chains have no incomparable pair of elements. 

Critical pairs can also be characterised as follows: $(u,v)$ is a critical pair of $\P$ if the addition of $u<v$ to the relations of $\P$ does not transitively force any other additional relation (or, equivalently, $v<u$ cannot be forced by adding any other relation). So $u<v$ could be seen as the canonical order of $u$ and $v$ in an LE $L$ of $\P$.  If $v < u$ holds in~$L$, we say that $L$ \emph{reverses} the critical pair $(u,v)$. If $L$ reverses some critical pair, it is \emph{reversing}. Here is the conjecture made in~\cite{felsner-led}:

\begin{conjecture}[Felsner, Reuter '99]
\label{conj}
Let $\P$ be a poset which is not a chain. Then in every diametral pair of LEs of $\P$, at least one of the two LEs is reversing.
\end{conjecture}

In Section~\ref{counter} we provide a counterexample to this conjecture. 
However, the conjecture turns out to be true in many special cases. We want to motivate this by looking at LE-graphs. We define a \emph{diametral vertex} as a vertex contained in a diametral pair of vertices. 

\begin{figure}[ht]
  	\begin{center}
     	\psfrag{1}{$1$}
    	\psfrag{2}{$2$}
    	\psfrag{3}{$3$}
    	\psfrag{4}{$4$}
    	\psfrag{a}{$x$}
    	\psfrag{b}{$y$}
    	\psfrag{P}{$\mathcal{M}$}
    	\psfrag{G}{$G(\mathcal{M})$}
    	\includegraphics[scale=0.4]{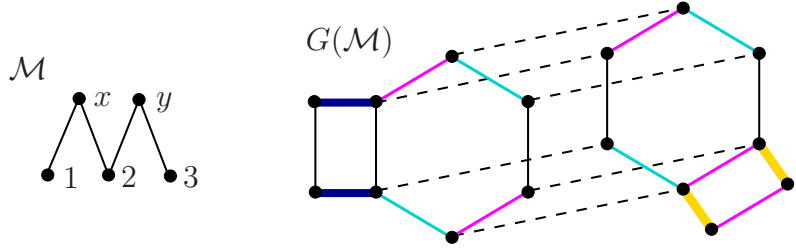}
	\caption{The poset $\mathcal{M}$ and its LE-graph: The 6-cycles in $G(\mathcal{M})$ correspond to the permutations of the antichain $1,2,3$. The dashed edges swap the pair $x,y$. The thick edges belong to the two critical pairs $(1,y)$ and $(3,x)$ of $\mathcal{M}$.}
	\label{LEG_M}
  	\end{center}
\end{figure}

Consider the LE-graph of $\mathcal{M}$ in Figure~\ref{LEG_M}, where we have highlighted the colours corresponding to critical pairs. 
 It looks as if these critical colours cut off the ``extremal parts'' of the graph. We can back up this intuition somewhat by observations about swap colours: If we remove the edges of one colour (say $ab$) from an LE-graph~$G(\P)$, then the graph falls into two components, the \emph{colour components}, corresponding to $\P+\{a<b\}$ and \mbox{$\P+\{a>b\}$}. The colour component $G(\P+\{a<b\})$ may be contained in some other colour component $G(\P+\{c<d\})$; this is the case exactly if $a<b$ forces $c<d$ (as for example $x<y$ forces $1<y$ in $\mathcal{M}$). Now for a critical pair $(u,v)$, we know that $v<u$ is not forced by any other pair, so $G(\P+\{v<u\})$ does not contain any other colour component. In this sense the critical colours are extremal colours, and it seems plausible that every diametral vertex should sit behind an extremal colour. This turns out not to be correct in all cases, but we will confirm it for several interesting special classes of posets. 

We make the following definition:

\begin{definition}
\label{diam_rev}
A poset $\P$ is \emph{diametrally reversing} if every diametral LE of~$\P$ is reversing.
\end{definition}

Note that chains are not diametrally reversing, because they contain no critical pair. This trivial observation will be crucial for the construction of the counterexample to Felsner and Reuter's conjecture.

In Section \ref{positive} we prove that several classes of posets are diametrally reversing, including posets of height $2$, interval posets and 3-layered posets. From the last class it follows that almost all posets are diametrally reversing.

\section{Complexity of Linear Extension Diameter}

For which cases is it easy to determine the LE-diameter of a given poset $\P$? If~$\P$ is a chain,
 then~$G(\P)$ consists of a single vertex, so its diameter is~$1$. If the poset is~$\mathcal{A}_n$, the antichain on
 $n$ elements, then a diametral pair consists of two permutations of the elements such that one is the reverse of
 the other. Hence $\textmd{led}(\mathcal{A}_n) =  {n \choose 2}$. 

If~$\P$ has dimension~$2$, then by definition there are two LEs such that every incomparable pair of elements appears in both orders, so $\textmd{led}(\P)$ equals the number of incomparable pairs $\textmd{inc}(\P)$ of~$\P$ and can thus be computed easily. Note that the width of a poset  is not greater than its dimension, thus determining the LE-diameter is also easy for all posets of width~$2$. 

In the following subsection we describe an algorithm that determines the LE-diameter of posets of width $3$, using a dynamic programming approach.

\subsection{Posets of Width $3$}

\begin{theorem}
	The LE-diameter of posets of width 3 can be computed in polynomial time.
\end{theorem}

\proof
Let a poset $\P$ of width 3 be given. By Dilworth's Theorem, it has a decomposition into three chains. Such a decomposition can be found in linear time with an algorithm of Felsner, Raghavan and Spinrad~\cite{felsner-dilworth}. Let $\P$ have $n$ elements, partitioned into the three chains $A, B$ and $C$. A convenient consequence is that $\P$ can have no more than $n^3$ downsets, because every downset is uniquely determined by an antichain of elements forming the maxima of the downset. Let $\P_D$ denote the poset induced by a downset~$D$ of $\P$. Our approach is to calculate the LE-diameter of $\P$ dynamically by calculating it for every $\P_D$ and re-using data in the process. 

Let us analyse what LEs of $\P$ can look like. In this proof we always read LEs from top to bottom, so the initial segment of an LE consists of its topmost vertices, to find the $i$-th element we count from the top, and so on. The first element of an LE $L$ can either be the top element of $A$, the top element of $B$ or the top element of $C$. Then there is an initial segment with only elements of that chain, until at a certain position the top element of a second chain appears. 

For chains $V,W,X,Y \in \{A,B,C\}$ and $i,j \in \mathbb{N}$ with $i,j \geq 2$, we define $\textmd{led}_D(VW, i, XY, j)$ as the maximum distance between two LEs $L_1, L_2$ of $\P_D$ such that: the maximal element of $L_1$ belongs to chain $V$, the second chain appearing is chain $W$, and the first element of chain $W$ appears at position~$i$ in $L_1$; the maximal element of $L_2$ belongs to chain $X$, the second chain appearing is chain $Y$, and the first element of chain $Y$ appears at position~$j$ in $L_2$. If there are no two such LEs, we set $\textmd{led}_D(VW, i, XY, j)= -\infty$.

Our plan now is to calculate the value defined above for every downset $D$ of $\P$, for every selection of the two first chains in $L_1$ and $L_2$, and for every position for the second chain to appear in $L_1$ and $L_2$. The number of values to compute is then bounded by $n^3 \cdot 6 \cdot 6 \cdot n \cdot n = O(n^5)$, thus polynomial. 

If $\P_D$ consists of only one chain or a chain plus one element,  then we can read off the desired value immediately. For all other cases we describe a recursive formula. So let us assume that $\P_D$ contains elements from all three chains, or that there are two chains from which $\P_D$ contains more than one element, and that we have computed all values for all downsets of $\P_D$. We want to compute $\textmd{led}_D(VW, i, XY, j)$ with $V,W,X,Y \in \{A,B,C\}$. Now we first choose a chain which appears in $VW$ as well as in $XY$. Since there are only three chains available, it must exist. Suppose it is chain $A$, then we will use the values for $D-a_1$ for our recursion, where $a_i$ is defined as the $i$-th element of chain $A$ (counted from the top) in $\P_D$. Define $D' :=D-a_1$. There are three cases to consider.

First we consider the case where $A$ is the chain appearing second in both $L_1$ and $L_2$ and the chains appearing first are different. We show a formula for $\textmd{led}_D(BA, i, CA, j)$; the formula for $\textmd{led}_D(CA, i,BA, j)$ is analogous. If $a_1  > b_{i-1}$ or $a_1  > c_{j-1}$,  then  $\textmd{led}_D(BA, i, CA, j) = - \infty$. Otherwise we know that the element~$a_1$ contributes exactly $i-1+j-1$ reversals to the distance between $L_1$ and $L_2$. Hence to calculate $\textmd{led}_D(BA, i, CA, j)$, we can maximise over all distances between two LEs $L_1', L_2'$ of $D'$ in which we can reinsert $a_1$ to get two LEs as specified. We call those LEs \emph{relevant}. What do our relevant LEs $L_1'$ and $L_2'$ look like? Of course, $L_1'$ starts with $i-1$ elements of chain $B$, and $L_2'$ with $j-1$ elements of chain $C$. So the first chain appearing in $L_1'$ and $L_2'$ is fixed. But after the initial segment, anything can happen in $L_1'$ and $L_2'$. Maybe the second appearing chain is again $A$ in both $L_1'$ and $L_2'$, in which case we are interested in the maximum distance over all possible positions of $a_2$, so let us set $\alpha =  \max_{r \geq i, p \geq j} \textmd{led}_{D'}(BA,r,CA,p)$. But maybe the distance can get bigger if the second chain in $L_1'$ is $C$, so we want to know $\beta =  \max_{r \geq i, p \geq j} \textmd{led}_{D'}(BC,r,CA,p)$; or if the second chain in $L_2'$ is $B$, in which case we are interested in $\gamma =  \max_{r \geq i, p \geq j} \textmd{led}_{D'}(BA,r,CB,p)$. The combination of these two cases yields the last case, where we look at $\delta =  \max_{r \geq i, p \geq j} \textmd{led}_{D'}(BC,r,CB,p)$. For $\textmd{led}_D(BA, i, CA, j)$ we then have to maximise over all these possibilities again, and we get

\[
	\textmd{led}_D(BA, i, CA, j)  =  i+j-2+ \max \{ \alpha, \beta, \gamma, \delta\}.
\]

In the second case $A$ appears first in one LE, say $L_1$, and second in the other, $L_2$. We will show a formula for $\textmd{led}_D(AB, i,CA, j)$. The formulae for $\textmd{led}_D(AC, i,BA, j)$ and also for $\textmd{led}_D(AB, i,BA, j)$ and $\textmd{led}_D(AC, i,CA, j)$ are built analogously. Again, if $b_1> a_{i-1}$ or $a_1 > c_{j-1}$ holds, we have $\textmd{led}_D(AB, i,CA, j) = - \infty$. Otherwise, we know that $a_1$ contributes $j-1$ reversals to the distance between $L_1$ and $L_2$. Now we need to distinguish the case that $i>2$ from the case $i=2$. If $i>2$ then we know that a relevant $L_1'$ starts with $i-2$ elements of $A$ and then lists $b_1$. So we only need to maximise over the possible second chains in $L_2'$, and obtain

\[
	\textmd{led}_D(AB, i, CA, j)  =  j-1+ \max \{ \alpha, \beta\},\quad \textmd{where}
\]
\[
	\alpha := \max_{p \geq j} \textmd{led}_{D'}(AB, i-1,CA,p) \quad \textmd{and} \quad
	\beta := \max_{p \geq j} \textmd{led}_{D'}(AB,i-1,CB,p).
\]

If $i=2$, then a relevant $L_1'$ starts with an element from $B$, and we need to maximise over all possible combinations of second chains in $L_1'$ and in $L_2'$:
\[
	\textmd{led}_D(AB, 2, CA, j)  =  j-1+ \max \{ \alpha, \beta, \gamma, \delta  \}, \quad \textmd{where}
\]
\begin{eqnarray*}
	\alpha &:=  \max\limits_{r \geq 2, p \geq j} \textmd{led}_{D'}(BA,r,CA,p), \quad 
	\beta &:= \max_{r \geq 2, p \geq j}\textmd{led}_{D'}(BA,r,CB,p), \\
	\gamma &:=  \max\limits_{r \geq 2, p \geq j} \textmd{led}_{D'}(BC,r,CA,p), \quad
	\delta &:= \max_{r \geq 2, p \geq j} \textmd{led}_{D'}(BC,r,CB,p).
\end{eqnarray*}

In the third case, $A$ is the first chain in both $L_1$ and $L_2$. We will show how to compute $\textmd{led}_D(AB, i, AC, j)$; the other cases for the second chains work analogously. The principle is the same, we only need to distinguish more cases now. First we observe that $a_1$ does not contribute any reversals here, but on the other hand we will always get a finite value for $\textmd{led}_D(AB, i, AC, j)$.  If both $i>2$ and $j>2$, we know exactly which chains come first and second in the relevant LEs; nothing changes but the position of the first element of the second chain: 
\[
	\textmd{led}_D(AB, i, AC, j) = \textmd{led}_{D'}(AB,i-1,AC,j-1).
\]
If one of the two parameters is exactly $2$, then we need to distinguish two cases as before. If $i=2$ and $j>2$, we have
\[
	\textmd{led}_D(AB, 2, AC, j)  = \max \{ \alpha, \beta\},\quad \textmd{where}
\]
\[
	\alpha := \max_{r \geq 2} \textmd{led}_{D'}(BA,r,AC,j-1)\quad \textmd{and} \quad
	\beta := \max_{r \geq 2} \textmd{led}_{D'}(BC,r,AC,j-1).
\]
 For $i>2$ and $j=2$, we obtain
\[
	\textmd{led}_D(AB, i, AC, 2)  = \max \{ \alpha, \beta\},\quad \textmd{where}
\]
\[
	\alpha :=\max_{p \geq 2} \textmd{led}_{D'}(AB,i-1,CA,p),  \quad \textmd{and} \quad
	\beta := \max_{p \geq 2} \textmd{led}_{D'}(AB,i-1,CB,p) .
\]
Finally, if both $i=2$ and $j=2$, we again have to consider four cases:
\[
	\textmd{led}_D(AB, 2, AC, 2)  = \max \{ \alpha, \beta, \gamma, \delta  \}, \quad \textmd{where}
\]
\begin{eqnarray*}
	\alpha &:= \max\limits_{r \geq 2, p \geq 2} \textmd{led}_{D'}(BA,r,CA,p), \quad 
	\beta &:= \max_{r \geq 2, p \geq 2} \textmd{led}_{D'}(BA,r,CB,p), \\
	\gamma &:= \max\limits_{r \geq 2, p \geq 2} \textmd{led}_{D'}(BC,r,CA,p), \quad
	\delta &:=\max_{r \geq 2, p \geq 2} \textmd{led}_{D'}(BC,r,CB,p).
\end{eqnarray*}

Hence we can compute the desired values for all downsets in time polynomial in $n$. The LE-diameter of the whole poset~$\P$ is now of course the maximum of the values for $\P_D=\P$ over all choices of the first and the second chain and all positions. \qed

Unfortunately we have found no way so far of generalising the above proof to an algorithm for determining the LE-diameter of posets of arbitrary fixed width. 

\subsection{General Posets}

In the following we will prove that it is NP-complete to determine the LE-diameter of a general poset. More precisely, we consider the following decision problem:

\vspace{.3cm}
\ni LINEAR EXTENSION DIAMETER \\
Input: Finite poset $\P$, natural number $k$ \\
Question: Are there two LEs of $\P$ with distance at least $k$?

\vspace{.3cm}

For the hardness proof we use a reduction of the following problem:

\vspace{.3cm}
\ni  BALANCED BIPARTITE INDEPENDENT SET\\
Input: Bipartite graph $G$, natural number $k$ \\
Question: Is there an independent set of size $2k$, consisting of $k$ vertices in each bipartition set?
\vspace{.3cm}

This problem is NP-complete as the equivalent problem BALANCED COMPLETE BIPARTITE SUBGRAPH is NP-complete \cite{gareyjohnson}. Before going into the proof, we need one more definition and a lemma. 

\begin{definition}
	In a poset $\P$, we call a subset $M$ of the elements a \emph{module} if they cannot be distinguished from the outside, that is, if for any $x \in \P \setminus M$ we have either $x > m \:\: \forall \, m \in M$ or $x < m \:\: \forall \, m \in M$ or $x ||m \:\: \forall \, m \in M$.
\end{definition}

The following lemma was already shown in \cite{felsner-led}; we provide a proof for completeness.
\begin{lemma}
	\label{modules}
For every poset $\P$ and module $M$ of $\P$, there is a diametral pair of LEs of $\P$ such that in both LEs of the diametral pair the elements of $M$ appear successively. 
\end{lemma}
\proof
Consider a diametral pair $L_1, L_2$ of LEs of $\P$. Each element $x \in \P$ contributes a certain number of reversals to the distance of $L_1$ and $L_2$, which is the number of elements $y \in \P$ such that the order of the pair $\{x,y\}$ is different in $L_1$ and $L_2$. Let $x$ be an element of $M$ which contributes a maximum number of reversals among all the elements of $M$. Now we move all the other elements of $M$ to the position of $x$ in $L_1$ and $L_2$ without changing their internal order. Since $M$ is a module, this cannot violate any relation of $\P$, thus the result will again be two LEs of~$\P$. By the choice of $x$, their distance cannot have decreased. So we have constructed a diametral pair such that the elements of $M$ appear consecutively in both LEs. \qed

\begin{theorem}
	LINEAR EXTENSION DIAMETER is NP-complete. 
\end{theorem}

\proof
The problem is in NP, because the distance of two given LEs can be checked in time quadratic in the number of elements of $\P$, and thus two LEs at distance $k$ are a certificate for a YES-instance. 

To show that the problem is NP-hard, suppose that an instance of BALANCED BIPARTITE INDEPENDENT SET is given: A bipartite graph $G=(A\cup B, E)$ and $k \in \mathbf{N}$. In a preprocessing step, we take a disjoint union of two copies $G_1=(A_1 \cup B_1, E_1)$ and $G_2=(A_2 \cup B_2, E_2)$ of $G$ and join all vertices of $A_1$ to all vertices of $B_2$ as well as all vertices of $A_2$ to all vertices of $B_1$. We call the resulting graph $G'=(A' \cup B', E')$. Then $G$ has a balanced independent set of size $2k$ exactly if $G'$ has two disjoint balanced independent sets of size $2k$. A further convenient fact is that, after this preprocessing step, we know that $k \leq |A'|/2, |B'|/2$. 

\begin{figure}[ht]
   \begin{center}
    \psfrag{a1}{$a_1$}
    \psfrag{b1}{$b_1$}
    \psfrag{a2}{$a_2$}
    \psfrag{b2}{$b_2$}
    \psfrag{ar}{$a_r$}
    \psfrag{bs}{$b_s$}
    \psfrag{A1}{$A_1$}
    \psfrag{A2}{$A_2$}
    \psfrag{Ar}{$A_r$}
   \psfrag{B1}{$B_1$}
    \psfrag{B2}{$B_2$}
    \psfrag{Bs}{$B_s$}
   \psfrag{C}{$C$}
    \psfrag{D}{$D$}
    \psfrag{G'}{$G'$}
   \psfrag{A'}{$A'$}
    \psfrag{B'}{$B'$}
    \psfrag{P}{$\P$}
    \psfrag{...}{$\ldots$}
    \includegraphics[scale=0.5]{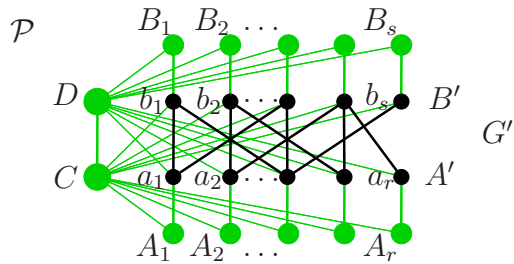}
    \caption{For the hardness proof, we build a poset from a bipartite graph.}
    \label{gadget}
   \end{center}
\end{figure}

Now we build a poset $\P$ starting from $G'$ by designating the vertices of $G'$ as \emph{black elements} of $\P$. The sets $A'$ and $B'$ each form an antichain in~$\P$, and an element of $B'$ is larger than an element of $A'$ exactly if they are adjacent in $G'$. Then we add the \emph{green elements} $A_1, A_2, \ldots, A_r$, $B_1, B_2, \ldots, B_s$ and $C$ and $D$ with relations as shown in Figure \ref{gadget}. 

Finally we form $\P^*$ from $\P$ by replacing the green elements by long chains. Let $n$ be the number of vertices of $G'$. Each $A_i$ and each $B_j$ is replaced by a chain of length $2n^4$, and $C$ and $D$ by chains of length $(2k-1)n^4$. All these chains form modules in $\P^*$, so we know by Lemma \ref{modules} that there is a diametral pair of LEs of $\P^*$ in which they appear successively. Since we are only interested in the distance between the LEs of a diametral pair, it suffices to consider such a diametral pair. This diametral pair corresponds to a pair of LEs of $\P$. The distance in $\P^*$ between the two LEs can be thought of as a \emph{weighted distance} between LEs of $\P$. From now on, we work in $\P$. 

We will analyse what a weighted diametral pair of LEs of $\P$ has to look like, and we will eventually see that its distance depends on the existence of a balanced independent set of $G'$. Recall that the distance between two LEs is the number of reversals between them. There are three types of pairs of elements of $\P$ that can be reversed: First, the pairs consisting of two black elements. Every such black/black reversal adds $1$ to the distance. We also call these reversals the \emph{unit reversals}. Second, there are the reversals of a black element with a green element. A black/green reversal contributes $2n^4$ or $(2k-1)n^4$ unit reversals to the distance, depending on the type of the green element involved. Finally there are reversals between two green elements, that is, between two $A_i$ or two $B_j$. Every such green/green reversal yields $4n^8$ unit reversals. 

We saw that the contribution of the three types of reversals differs a lot. In fact, the total gain of all possible reversals of one type yields still less than one reversal of the next bigger type: There are $n$ black elements, so the number of black/black reversals is at most ${n \choose 2}$. This is $\theta(n^2)$, far less than the $\theta(n^4)$ that one black/green reversal yields. And how many of those can there be in total? Even if the black elements were incomparable to the whole rest of the poset, the black/green reversals could altogether contribute only $\theta(n^6)$ unit reversals, again far less than the $\theta(n^8)$ of a single green/green reversal. 

The consequence of this is that we can analyse the LEs forming a weighted diametral pair of $\P$ in three separate steps. First we check how the green elements have to be ordered in the LEs to yield a maximum number of green/green reversals. In the second step, we fill in the black elements in order to get a maximum number of black/green elements, knowing that these cannot influence the order of the green elements, because they cannot contribute enough. In the last step, we use the remaining freedom to order the black elements so that they add a maximum number of unit reversals. Note that the first two steps do not depend on the particular graph $G'$. We define the \emph{base distance} to be the weighted distance we can achieve between two LEs of $\P$ independently of the particular given graph $G'$; put differently, it is the LE-diameter of $\P^*$ built from a complete bipartite graph $G'$. In the last step we will see that the existence of a balanced independent set determines how much $\textmd{led}(\P^*)$ exceeds the base distance. 

\begin{figure}[ht]
   \begin{center}
    \psfrag{a1}{$a_1$}
    \psfrag{b1}{$b_1$}
    \psfrag{a2}{$a_2$}
    \psfrag{b2}{$b_2$}
    \psfrag{ar}{$a_r$}
    \psfrag{bs}{$b_s$}
    \psfrag{A1}{$A_1$}
    \psfrag{A2}{$A_2$}
    \psfrag{Ar}{$A_r$}
   \psfrag{B1}{$B_1$}
    \psfrag{B2}{$B_2$}
    \psfrag{Bs}{$B_s$}
   \psfrag{C}{$C$}
    \psfrag{D}{$D$}
    \psfrag{G'}{$G'$}
   \psfrag{A'}{$A'$}
    \psfrag{B'}{$B'$}
    \psfrag{P}{$\P'$}
    \psfrag{...}{$\ldots$}
    \includegraphics[scale=0.5]{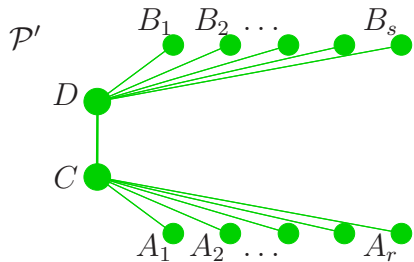}
    \caption{In the first step, we only consider the green part $\P'$ of $\P$.}
    \label{greenpart}
   \end{center}
\end{figure}

We start with the first step. For this, it is enough to look at the poset~$\P'$ induced only by the green elements, so $\P' = \P \setminus G'$. It consists of an antichain of minima, formed by the $A_i$, the two elements $C$ and $D$ which are comparable to all other elements of $\P'$, and an antichain of maxima, formed by the $B_j$; see Figure~\ref{greenpart}. 

The poset $\P'$ is two-dimensional, and the diametral pairs consist of two LEs $L_1'$, $L_2'$ such that the $A_i$ and $B_j$ come in an arbitrary order in $L_1'$, and in the opposite order in $L_2'$. Thus up to labelling the $A_i$ and $B_j$, a diametral pair of LEs of $\P'$ has the form
\begin{eqnarray*}
	L_1' &=&A_1 A_2 \ldots A_r C D B_1 B_2 \ldots B_s \\
	L_2 '&=&A_r A_{r-1} \ldots A_1 C D B_s B_{s-1} \ldots B_1.
\end{eqnarray*}
Hence the green/green reversals in a diametral pair contribute a total number of $\left({r \choose 2}+{s \choose 2}\right) \cdot 4n^8$ unit reversals to the base distance. 

\begin{figure}[ht]
   \begin{center}
    \psfrag{a1}{$a_1$}
    \psfrag{b1}{$b_1$}
    \psfrag{a2}{$a_2$}
    \psfrag{b2}{$b_2$}
    \psfrag{ar}{$a_r$}
    \psfrag{bs}{$b_s$}
    \psfrag{A1}{$A_1$}
    \psfrag{A2}{$A_2$}
    \psfrag{Ar}{$A_r$}
   \psfrag{B1}{$B_1$}
    \psfrag{B2}{$B_2$}
    \psfrag{Bs}{$B_s$}
   \psfrag{C}{$C$}
    \psfrag{D}{$D$}
    \psfrag{G'}{$G'$}
   \psfrag{A'}{$A'$}
    \psfrag{B'}{$B'$}
    \psfrag{P}{$\P'$}
    \psfrag{...}{$\ldots$}
    \includegraphics[scale=0.5]{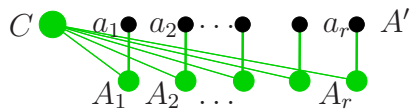}
    \caption{In the second step, we look at the lower half $\bar{\P}$ of $\P$.}
    \label{lowerhalf}
   \end{center}
\end{figure}

For the second step we need to take the elements of $G'$ into account. We concentrate on the \emph{lower half} $\bar{\P}$ of $\P$, induced by the elements of $A_1 \cup \ldots \cup A_r \cup C \cup a_1 \cup \ldots a_r$ (see Figure~\ref{lowerhalf}). 

The subposet induced by the remaining elements is the \emph{upper half}. The only incomparabilities between the upper and the lower half occur between two black elements, and we are not interested in black/black reversals in the second step. Otherwise the lower half and the upper half are symmetric, so the number of black/green reversals we can achieve in total is twice the number of possible black/green reversals in the lower half. 

Let us look at the lower half. Each $a_i$ is incomparable to $C$ and all $A_j$ with $j\neq i$. One option to insert $a_i$ into $L_1'$ and $L_2'$ yielding many reversals is to place it right above its predecessor $A_i$ in both LEs, thus reversing with all $A_j$, $j \neq i$, but not with $C$. The other option is to insert $a_i$ above $C$ (and thus also above all $A_j$) in one LE, and place it as low as possible in the other. Note for this option that $a_i$ still has to be above $A_i$, and thus it cannot be reversed with the $A_j$ which are below $A_i$ in the other LE. 

In the first option, $a_i$ contributes $(r-1) \cdot 2n^4$ unit reversals. In the second option, we win $(2k-1)n^4$ unit reversals from $C$, but lose $2n^4\cdot \min\{ i-1, r-i  \} $ unit reversals from the $A_j$. Recall for the following that $k \leq r/2$. If $i \leq k$, then the minimum is attained by $i-1$, and $2n^4\cdot (i-1)  < (2k-1)n^4$. Thus it is best to choose the second option and insert $a_i$ above $C$ in $L_2'$. We say that the elements between $C$ and $D$ are \emph{in the middle} of an LE of $\P$. If $i \geq r-k+1$, then $r-i$ attains the minimum, and again $2n^4\cdot (r-i) < (2k-1)n^4$. In this case it is best to insert $a_i$ in the middle of $L_1'$. Conversely, if $k+1 \leq i \leq r-k$, then $2n^4\cdot \min\{ i-1, r-i  \} > (2k-1)n^4$, thus the first option is best. Hence exactly $k$ of the $a_i$ will be inserted in the middle of $L_1'$ and $L_2'$, respectively.

The analysis for the upper half can be made in exactly the same way, so $k$ of the $b_j$ end up in the middle of the $L_1'$ and $L_2'$, respectively. Up to reversals in the middle, we can now write down what a diametral pair $L_1, L_2$ of LEs of the whole poset $\P$ looks like. 
\begin{eqnarray*}
	L_1 & = & A_1 a_1 A_2 a_2 \ldots A_r C a_{r-k+1} \ldots a_{r-1} a_r  b_1 b_2 \ldots b_k  D B_1 B_2 \ldots  b_{s-1} B_{s-1} b_s B_s \\
	L_2 & = & A_r a_r A_{r-1} a_{r-1} \ldots A_1 C a_k \ldots a_2 a_1  b_s b_{s-1} \ldots b_{s-k+1}  D B_s B_{s-1}\ldots b_2 B_2 b_1 B_1
\end{eqnarray*}
Note that the chosen order of the $A_i$ and $B_j$ determines which elements $a_i$ and $b_j$ end up in the middle of $L_1$ and $L_2$.

The black/green reversals contributed by the lower half amount to the following number of unit reversals: 
\begin{eqnarray*}
	 r(r-1) \cdot 2n^4   &  - 2\left( \sum_{i=1}^k (i-1)\right)\cdot2n^4 & + 2k \cdot (2k-1)n^4  \\
	=  r(r-1) \cdot 2n^4 & - k(k+1)\cdot2n^4 & + 2k \cdot (2k-1)n^4.
\end{eqnarray*}
For the number of contributed black/green reversals of the upper half, we only have to replace~$r$ by~$s$. 

The third step is now easy. To analyse how many black/black reversals we can obtain, we first observe that it is clearly possible to completely reverse the order of the $a_i$ as well as the order of the $b_j$. This adds another ${r \choose 2} + {s \choose 2}$ to the base distance.

Now we put our calculations together to find the base distance $d$ between two diametral LEs of $\P$. 
\begin{eqnarray*}
d = &\sum_{t =r,s}  {t \choose 2}\cdot 4n^8  + t(t-1) \cdot 2n^4  - k(k+1)\cdot2n^4  + 2k \cdot (2k-1)n^4 + {t \choose 2} \\
=  &\left( {r \choose 2} + {s \choose 2 } \right) \cdot \left( 2n^4 +1  \right)^2 - 2k(k+1)\cdot2n^4  + 4k \cdot (2k-1)n^4
\end{eqnarray*}

For the only other possible black/black reversals, we have to check if some of the $a_i$ can be brought above some of the $b_j$. This can only happen between the $a_i$ and $b_j$ in the middle of $L_1$ and $L_2$. So there are at most  $2k^2$ additional unit reversals to be won. This number can be obtained exactly if $G'$ contains two disjoint balanced independent sets of size $2k$.  Because if it does, then we can choose an order of the $A_i$ and $B_j$ so that the vertices of one independent set end up in the middle of $L_1$, and the vertices of the other in the middle of $L_2$, and then bring all the involved $a_i$ above all the involved $b_j$. On the other hand, if we can win these additional $2k^2$ unit reversals, then clearly there have to be two disjoint balanced independent sets in $G'$.

Recall that $G'$ has two disjoint balanced independent sets of size $2k$ exactly if $G$ has one. So we conclude that~$G$ has a balanced independent set of size~$2k$ exactly if~$\P^*$ has two LEs of distance at least~$d+2k^2$. Since it is clearly possible to build $\P^*$ in time polynomial in the size of $G$, this proves the hardness of LINEAR EXTENSION DIAMETER. \qed

\section{The Counterexample}
\label{counter}

In this section we construct a counterexample to Conjecture~\ref{conj}. We will first give a nontrivial example of a poset which is not diametrally reversing, because it provides a simpler version of the construction. The idea in both examples is to replace some elements of a Boolean lattice by long chain modules. The chains function like a weight on the elements they replace, and we can therefore use them to manipulate the behaviour of diametral pairs. At the same time, the chains do not add any new critical pairs, as the following lemma shows. 

\begin{lemma}
\label{crit_bij}
	Let $\P$ be a finite poset, and let $\P'$ arise from it by replacing each element $x$ of $\P$ by the chain $X= x_1\leq x_2 \leq \ldots \leq x_{k(x)}$, where all of these elements might be equal.  Then by mapping a pair $(x,y)$ of elements of $\P$ to the pair $(x_1, y_{k(y)})$ of elements of $\P'$, we obtain a bijection between the critical pairs of $\P$ and the critical pairs of $\P'$.
\end{lemma}

\proof
Observe that, since all the introduced chains are modules, we have $v<w$ in $\P$ exactly if all elements of $V$ are smaller than all elements of $W$ in $\P'$, and $v||w$ in $\P$ exactly if all elements of $V$ are incomparable to all elements of $W$ in $\P'$. 

If $(x,y)$ is a critical pair of $\P$, then $x||y$ in $\P$ and hence $x_1 ||y_{k(y)}$ in $\P'$. From $\downarrow_{\P} \!\! x \subseteq \downarrow_{\P} \!\! y$ we deduce $\downarrow_{\P'} \!\! x_1 \subseteq \downarrow_{\P'} \!\! y_{k(y)}$. Analogously $\uparrow_{\P} \!\! y \subseteq \uparrow_{\P} \!\! x$ yields $\uparrow_{\P'} \!\! y_{k(y)} \subseteq \uparrow_{\P'} \!\ \!x_1$, and thus  $(x_1, y_{k(y)})$ is a critical pair of $\P'$. 

Now let $(x_i, y_j)$ be a critical pair of $\P'$. Then we claim that $i=1$ and $j=k(y)$ must hold. Since $x_i || y_j$, all pairs of elements from the two chains are incomparable. Therefore each element in $X$ other than $x_1$ has a predecessor which is not a predecessor of $y_j$, namely, $x_1$. In the same way, each element in $Y$ other than $y_{k(y)}$ has a successor which is not a successor of $x_i$, namely, $y_{k(y)}$. Hence all critical pairs of $\P'$ have the form $(x_1, y_{k(y)})$. On the other hand, if $(x_1, y_{k(y)})$ is critical in~$\P'$, then clearly $(x,y)$ is critical in~$\P$, and thus the defined mapping is a bijection. \qed

\begin{figure}[h]
  \begin{center}
     \psfrag{0}{$\emptyset$}
    \psfrag{1}{$1$}
    \psfrag{2}{$2$}
    \psfrag{3}{$3$}
    \psfrag{4}{$4$}
     \psfrag{12}{$12$}
    \psfrag{13}{$13$}
   \psfrag{14}{$14$}
    \psfrag{23}{$23$}
   \psfrag{24}{$24$}
    \psfrag{34}{$34$}
    \psfrag{123}{$123$}
    \psfrag{134}{$134$}
    \psfrag{234}{$234$}
    \psfrag{124}{$124$}
   \psfrag{1234}{$1234$}
	\includegraphics[scale=0.45]{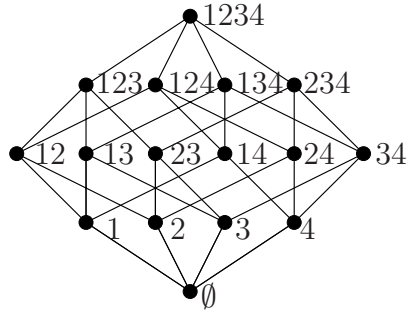} 
	\caption{The Boolean lattice $B_4$.}
	\label{b4}
\end{center}
\end{figure}

We will also need some facts about the Boolean lattice $B_n$, which is the poset on the subsets of $[n]=\{1,\ldots, n\}$, ordered by inclusion. Figure~\ref{b4} shows the Hasse diagram of $B_4$. It is a well-known fact that the critical pairs of~$B_n$ are formed by the atoms with the corresponding coatoms. We prove a generalisation of this fact which we will need later.

\begin{lemma}
\label{Bn_crit}
	Let $\P$ be a subposet of $B_n$, induced by a set of subsets of $[n]$ which contains all atoms $i$ and all coatoms $[n]-i$, where $i=1,\ldots, n$. Then the critical pairs of $\P$ are exactly the pairs $(i, [n]-i)$. 
\end{lemma}
\proof
Suppose two subsets $S$ and $T$ of $[n]$ form a critical pair. Then they are incomparable, so $S \not\subseteq T$. Therefore $S$ contains an atom $i$ which is not contained in $T$. Now if $S \neq i$, then $i$ would be a predecessor of $S$ which is no predecessor of $T$, contradicting that $(S,T)$ is a critical pair. Hence $S=i$. Since $i \notin T$, we know that $T \subseteq [n]-i$. But if $T$ is a proper subset of the coatom $[n]-i$, then this coatom forms a successor of $T$ which is not a successor of $S$, a contradiction. Thus we have $T = [n]-i$, and the critical pairs of $\P$ are as claimed. \qed

\begin{theorem}
	Let $B_4^*$ be the poset resulting from $B_4$ if the two-element sets, or \emph{doubles}, are replaced with chains of length 3. Then $B_4^*$ is not diametrally reversing. 
\end{theorem}

\proof 
Let us first think of $B_4^*$ as the poset obtained from $B_4$ by replacing the doubles with long chains. Let $w$ be their length; later we will show that $w=3$ suffices. 

Note that the introduced chains all form modules in $B_4^*$. We are interested in pairs of LEs of $B_4^*$ with large distance. Since we can always make the chain modules appear successively in a pair of LEs without lowering the distance (cf.~Lemma~\ref{modules}), we will first look at such LEs. 
As in the hardness proof, the distance between two LEs of $B_4^*$ can be thought of as a \emph{weighted distance} between LEs of~$B_4$. Reversing two doubles yields $w^2$ \emph{unit reversals}. If $w$ is chosen large enough, then reversing as many doubles as possible has priority over all other reversals. Thus in a diametral pair $L_1, L_2$ of $B_4^*$, the doubles appear in some order in $L_1$ and in the opposite order in $L_2$. 

Let $L$ be a reversing LE of $B_4^*$. By Lemmas~\ref{crit_bij} and \ref{Bn_crit} we know that the critical pairs of $B_4^*$ are exactly the pairs $(i, [4]-i)$ for $i=1,2,3,4$. Assume without restriction that $234<1$ in~$L$.  Then we know that $\{23,24,34\} <234<1< \{12,13,14\}$ in $L$. Let $L'$ be an LE of $B_4^*$ which has maximum distance to $L$. For large $w$ we have $ \{12,13,14\} < \{23,24,34\}$ in $L'$. Suppose for contradiction that $L'$ reverses the critical pair $(j, [4]-j)$. Then the three doubles larger than~$j$ must appear in $L'$ after the three doubles smaller than $[4]-j$. But the three last doubles in $L'$ do not have an atom in common. Hence $L'$ cannot be reversing. We conclude that if a diametral LE of $B_4^*$ is reversing, then no diametral partner of it is. 

With a closer analysis we can bound $w$. Here are two LEs of $B_4$ which reverse all pairs of doubles: 
\begin{eqnarray*}
	L_1 & = & \emptyset \,\, 1 \,\,2  \,\,12  \,\,3  \,\,13  \,\,23  \,\,123  \,\,4 \,\,14  \,\,24  \,\,124 \,\,34 \,\,134 \,\,234 \,\, 1234 \\
	L_2 & = & \emptyset  \,\,4 \,\, 3 \,\, 34  \,\, 2 \,\,24  \,\,1  \,\, 14 \,\,23  \,\, 234 \,\,13 \,\,134 \,\,12 \,\, 124  \,\, 123  \,\,1234
\end{eqnarray*}
Their weighted distance is ${6 \choose 2}w^2 + 14w +13$. If we consider two reversing LEs, then as shown above their quadratic term is smaller. For the linear term, we observe that the largest atom in an LE of $B_4$ can be larger than at most three doubles, and the second largest atom larger than at most one double. Analogously the smallest coatom can be smaller than at most three doubles, the second smallest coatom smaller than at most one. Therefore the linear term is at most $16w$. The constant term is bounded by 14 by counting six reversals among atoms and among coatoms, respectively, and two reversals for the two critical pairs. We just have to choose $w$ such that ${6 \choose 2}w^2 + 14w +13 > \left({6 \choose 2}-1 \right) w^2 + 16w+14$, and this holds for $w=3$. 

So there is an actual gap between the distance of two reversing LEs of~$B_4^*$ and the distance of a diametral pair. This gap remains even if we do not insist that the chain modules appear successively, so we conclude that $B_4^*$ is not diametrally reversing.
\qed

Note that $B_4^*$ is a graded poset (cf.~Definition~\ref{graded}), so not all graded posets are diametrally reversing. Now we refine the ideas of the above construction to disprove Conjecture~\ref{conj}. 

\begin{theorem}
	There is a poset $\P^*$ such that no LE of $\P$ contained in a diametral pair is reversing. 
\end{theorem}

\begin{figure}[ht]

  \begin{center}
    \includegraphics[scale=0.4]{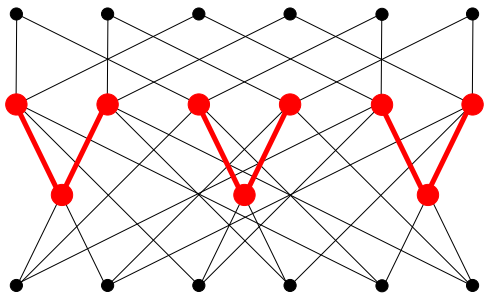}
   \caption{The counterexample to Felsner and Reuter's conjecture: the fat elements represent long chains.}
  \label{counterex}
	\end{center}
\end{figure}

\proof
We consider the subposet~$\P$ of the six-dimensional Boolean lattice~$B_6$ induced by all the atoms, the three doubles $12,34,56$, no triples, the six quadruples $1235, 1246, 1345, 2346, 1356, 2456$, and all the coatoms, see Figure~\ref{counterex}. We replace the doubles and the quadruples in~$\P$ by chains of length~$w$, which we will specify later. The resulting poset is our counterexample~$\P^*$. 

Again we are interested in pairs of LEs of $\P^*$ with large distance, and we first look only at LEs in which the chain modules appear successively (cf.~Lemma~\ref{modules}). We work with $\P$ instead of $\P^*$, treating the chains of length~$w$ as single elements, called \emph{red elements}, which have weight~$w$. We make~$w$ so big that gaining a maximum number of those red-red reversals, each contributing $w^2$ unit reversals, is more desirable than anything else. Let us call the poset induced by the red elements~$\P_{\textmd{red}}$. A maximum number of red-red reversals is achieved by taking two LEs of $\P^*$ which restricted to the red elements form a diametral pair of~$\P_{\textmd{red}}$. Now $\P_{\textmd{red}}$ is two-dimensional, so $\textmd{led}(\P_{\textmd{red}}) = \textmd{inc}(\P_{\textmd{red}})={3 \choose 2} + {6 \choose 2} + 3\cdot4 = 30$. Since we can extend an LE of $\P_{\textmd{red}}$ to an LE of $\P^*$, we obtain $\textmd{led}(\P^*)>30w^2$. 

Which distance can a pair of LEs achieve in which at least one of them is reversing? Again by Lemmas \ref{crit_bij} and \ref{Bn_crit} we know that the critical pairs of $\P^*$ are exactly the pairs $(i, [6]-i$), $i=1,\ldots, 6$. So let $(L_1, L_2)$ be a pair of LEs of $\P^*$, and let $L_1$ reverse the critical pair $(6, 12345)$, say. Then for the red elements we have $\{12, 34, 1235, 1345 \} < \{56, 1246, 2346, 1356, 2456\}$ in $L_1$. If there is an LE of~$\P_{\textmd{red}}$ contained in a realiser which respects these relations, then its partner in the realiser has to fulfill $1246<34$ and $2346<12$ at the same time, which is impossible. Thus there are no more than 29 red-red reversals between $L_1$ and $L_2$. 

Let us bound the number of other reversals between $L_1$ and $L_2$. Each atom of $\P^*$ is incomparable to four red elements, and each coatom to five. If we add all possible reversal of atoms and coatoms among themselves and finally the reversed critical pair, we see that the distance between $L_1$ and $L_2$ is at most $29w^2+6\cdot4w+6\cdot5w+6!+6!+1=29w^2+54w+1441$.

We choose $w$ so big that $30w^2>29w^2+54w+1441$, e.g., $w=100$. Then there is a gap between the distance we can achieve when reversing a critical pair and the distance of a diametral pair. Again we cannot gain anything by relaxing the condition that the chain modules have to appear successively, so we conclude that every diametral pair of $\P^*$ consists of two non-reversing LEs. \qed

Note that our construction disproves the conjecture in a very strong sense: Not only have we shown that not every diametral pair of our example contains a reversing LE, but in fact no diametral pair at all does.

\section{Diametrally Reversing Posets}
\label{positive}

In this section we will present a number of classes of diametrally reversing posets. Consequently all these posets fulfill Felsner and Reuter's Conjecture. We start with the well-understood class of 2-dimensional posets.

\begin{proposition}
\label{2dim}
Every poset of dimension $2$ is diametrally reversing. 
\end{proposition}
\proof
If $\P$ is a 2-dimensional poset, then the diametral pairs of LEs of $\P$ are exactly the realisers of $\P$. Now it is a classical result that a collection $L_1, L_2, \ldots, L_k$ of LEs forms a realiser if and only if every critical pair of $\P$ is reversed in at least one of them (cf.~\cite{trotter92}). If a non-reversing diametral LE $L$ would exist, then its partner $L'$ in a diametral pair would have to reverse all critical pairs of $\P$. But in this case, $L'$ would be a realiser of $\P$ all by itself. This is a contradiction to $\P$ being 2-dimensional.
\qed

\subsection{Modules}
\label{twinsnmodules}

What can we say about posets which are in some way constructed from smaller posets? One possibility is that $\P$ is a series composition of smaller posets $\P_i$. Then every LE of $\P$ is a concatenation of LEs of the $\P_i$, and the LE-graph of $\P$ is the cartesian product of the $G(\P_i)$. Therefore it is easy to see that it suffices that one of the $\P_i$ is diametrally reversing to make the whole poset $\P$ diametrally reversing. It turns out that this is even true for general modules:

\begin{proposition}
\label{M_suffices}
Let $\P$ be a poset containing a module $M$. If $M$ is diametrally reversing, then so is $\P$. 
\end{proposition}
\proof 
By Lemma~\ref{modules}, there is a diametral pair of LEs of $\P$ in which the elements of $M$ appear successively. Restricting the two LEs of such a pair to $M$ clearly yields a diametral pair of LEs of $M$. But even if in a diametral pair $L_1,L_2$ of LEs of $\P$ the elements of $M$ do not appear successively, they must contribute the same number of reversals, and hence the restriction of $L_1$ and $L_2$ to $M$ again yields a diametral pair of LEs of $M$. Now since $M$ is a module, the critical pairs of $M$ stay critical in $\P$. Hence if every diametral LE of $M$ reverses a critical pair, then also every diametral LE of $\P$ reverses a critical pair. \qed

A special case of a module is a \emph{twin}: two elements with the same downset and the same upset. Note that the two elements in a twin must be incomparable. They form a critical pair in both orders. So any LE reverses one of these two critical pairs, and we have the following very useful result:

\begin{corollary}
\label{twins}
Every poset containing a twin is diametrally reversing. 
\end{corollary}

\subsection{Interval Orders}

In this section we prove that interval orders are diametrally reversing.

\begin{definition}
A poset is an \emph{interval order} if its elements can be represented by intervals on the real line such that $u< v$ if and only if the interval representing $u$ is completely left of the interval representing $v$. If this can be done while all intervals have length $1$, we speak of a \emph{unit interval order}.	
\end{definition}

Note that the critical pairs of interval orders correspond to pairs of intervals which do not contain each other, ordered from left to right. By slight abuse of notation, we will not differentiate between an element and the corresponding interval in the following. We will first prove that unit interval orders are diametrally reversing. We need this result for the proof of the general case. The proof of the special case also contains basic ideas we will use repeatedly in later proofs.

\begin{proposition}
\label{unit}
Every unit interval order which is not a chain is diametrally reversing.
\end{proposition}
\proof
Let $\P$ be a unit interval order. We can assume that $\P$ is not an antichain, otherwise we are done by Proposition~\ref{2dim}. Consider a unit interval representation of $\P$. Note that the left-endpoint-order of the intervals is the same as the right-endpoint-order, since all intervals have the same length. As a consequence, each incomparable pair $u,v$ is a critical pair $(u,v)$, where $u$ is represented by the interval more to the left. So in fact there is only one non-reversing LE: the left-endpoint-order $L$. We will show that $L$ cannot be diametral. 

Cover the plane with a grid of vertical lines $a, b, c, \ldots$ of distance 1. Let us assume all intervals are closed on the right end and open on the left end. Then every interval intersects exactly one grid line. We place the grid in such a way that $a$ hits the right endpoint of the leftmost interval. Let us denote the elements intersecting line $a$ with $a_1, a_2, \ldots, a_k$,  ordered as in $L$, the elements intersecting line $b$ with $b_1, b_2, \ldots$, and so on. 
Note that we have $a_1 < b_1$ in $\P$, because if they would be incomparable, $b_1$ had to intersect line $a$, a contradiction. Also observe that if we would have $a_k < b_1$ in $\P$, then the antichain of the $a_i$ would form a (serial) module in $\P$, containing lots of twins. So by Corollary~\ref{twins} we can assume that $a_k || b_1$. 

Now assume for contradiction that $L$ forms a diametral pair with $L'$. In~$L$, element $a_k$ is adjacent to $b_1$, hence we must have $b_1 < a_k$ in $L'$, otherwise we could increase the distance of the two LEs by exchanging $a_k$ and $b_1$ in $L$. 

Similarly, the order of the $a_i$ can be chosen freely in $L$, since they form an antichain of successive elements. Therefore the order of the $a_i$ in $L$ must be the reverse of their order in $L'$, otherwise we could construct a pair of LEs with larger distance. 

So we conclude that in $L'$ we must have $b_1< a_k$ and $a_k < a_1$, which implies $b_1 < a_1$. But this is a contradiction because $a_1 < b_1$ in $\P$. 
\qed 

In the general case, we will use the same idea: We try to show that if, in a pair of LEs, one LE is not reversing, then we can always construct two LEs with larger distance. 

For the proof we want to use a canonical representation of interval orders. We need the following characterisation: A poset~$\P$ is an interval order if and only if the maximal antichains of $\P$ can be linearly ordered such that, for each element $v \in \P$, the maximal antichains containing $v$ occur consecutively (see e.g.~\cite{moehring-tractable}). The ``only if'' direction is immediately clear when considering an interval representation. For the other direction, we take a linear order of the antichains as in the characterisation and define $\ell(v)$ and $r(v)$ as the indices of the first and last antichains containing $v$. Then $(\ell(v), r(v))_{v \in \P}$ defines a representation of $\P$ by open intervals with integer endpoints which will be our \emph{canonical interval representation}.

\begin{figure}[ht]
   \begin{center}
    \psfrag{a}{$a$}
    \psfrag{b}{$c$}
    \psfrag{c}{$b$}
    \psfrag{d}{$d$}
    \psfrag{p}{$e$}
    \psfrag{s}{$w$}
    \psfrag{x}{$x$}
    \psfrag{y}{$y$}
    \psfrag{z}{$z$}
    \includegraphics[scale=0.5]{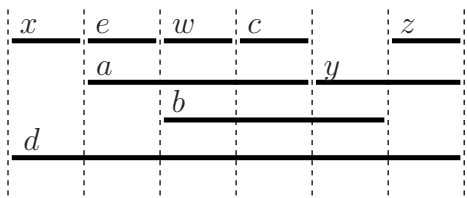}
    \caption{Example of an interval order in canonical representation.}
    \label{standard-rep}
   \end{center}
\end{figure}

Figure~\ref{standard-rep} shows an example where we have drawn vertical lines to mark the integers which constitute left or right endpoints of intervals. An observation important for the proof is that no intervals start or end between the lines, but at every line except for the last one an interval starts, and at every line except for the first one an interval ends. 

Now we have enough tools at hand to go into the proof. 

\begin{theorem}
Every interval order which is not a chain is diametrally reversing.
\end{theorem}
\proof
Let $\P$ be an interval order given in canonical representation. Let $L$ be a diametral LE of $\P$ with diametral partner $L'$. Suppose for contradiction that $L$ is non-reversing, that is, every critical pair $(x,y)$ of~$\P$ appears in the canonical order $x<y$ in $L$. A non-reversing LE of the example in Figure~\ref{standard-rep} is given by $L=xewabcdyz$. We will construct a pair of LEs with larger distance, contradicting the fact that $L,L'$ is a diametral pair. 

If there are no two intervals in $\P$ such that one contains the other, then the left-endpoint-order of the intervals equals the right-endpoint-order, so $\P$ is a unit interval order and we are done by the previous proposition. Thus there are pairs of intervals containing each other; these are exactly the pairs which are not critical.

Let $c$ be the first element in $L$ which appears after an element $a$ with $r(a) \geq r(c)$. If there is no such $c$, then we apply the proof to $L$ read backwards, exchanging left and right endpoints. Note that since $L$ is non-reversing, it follows that $\ell(a)< \ell(c)$. Let $b$ be the element appearing immediately before $c$ in $L$. Our plan is to move $b$ down in $L$ and up in $L'$ and thus create a pair of LEs with larger distance. 

We want to place $b$ \emph{low} in $L$, that is, into a position such that all the elements before $b$ in $L$ are predecessors of $b$. To show that this is possible we claim that the predecessors of $b$ form an initial segment of $L$. Indeed, if we consider an element $d$ incomparable to $b$, then clearly it has a larger right endpoint than any predecessor $e$ of $b$. Hence if $d$ would appear before $e$ in~$L$, we would have chosen $e$ as $c$, a contradiction. Therefore we can safely move~$b$ down in $L$ and place it right after its first predecessor, obtaining an LE $\bar{L}$. Let $S$ be the set of elements that $b$ passes, that is, that are smaller than $b$ in $L$ and larger than $b$ in $\bar{L}$. 

Now we move $b$ up in $L'$, that is, we move it just before its lowest successor in $L'$. We call the resulting LE $\bar{L}'$. First we show that the distance between $\bar{L}$ and $\bar{L}'$ is not smaller than the distance between $L$ and $L'$. Note that $b$ is low in $\bar{L}$ and thus any element that $c$ passes when moving up in $L'$ only increases the distance between the two LEs. If we can ensure that in $\bar{L}'$, the element $b$ is larger than all elements of $S$, then we have accounted for the elements that $b$ passed in $L$. But this can be done: By the choice of $c$, all elements which are smaller than $b$ in $L$, in particular the elements of $S$, have a smaller right endpoint than $b$. Therefore every successor of $c$ is also a successor of all elements in $S$. Hence if we move $b$ up in $L'$ and place it just before its lowest successor, we will have passed all elements of $S$. 

It remains to show that the distance between $\bar{L}$ and $\bar{L}'$ is actually larger than the distance between $L$ and $L'$. We will identify a pair of elements that is reversed between the first pair of LEs, but not between the latter. First let us concentrate on the pair $b,c$. 
By the choice of $c$, we must have $r(b) \geq  r(a) \geq r(c)$. It follows that $b$ and $c$ are incomparable. If $\ell(b) \geq \ell(c)$, then  $(c,b)$ is a critical pair which is reversed in $L$, a contradiction. Thus we have $\ell(b) < \ell(c)$.

Now from the properties of the canonical interval representation it follows that there exists an element $w$ which is incomparable to $b$ but smaller than~$c$. What is the order of $b,c$ and $w$ in $L$ and $L'$? In $L$ we know that $b$ and $c$ appear adjacently, so $w$ has to appear before both of them. In $L'$ the order of $b$ and $c$ must be reversed, otherwise we could construct a pair of LEs with larger distance immediately by changing their order in $L$. Hence in $L'$ we have $w<c<b$, and thus the pair $w,b$ is not reversed between $L$ and~$L'$. But $b$ is low in $\bar{L}$, and therefore $b<w$ in $\bar{L}$. Therefore $w,b$ is reversed between $\bar{L}$ and $\bar{L}'$ and we have found the desired pair. This shows that the non-reversing LE $L$ cannot be diametral. \qed

\subsection{3-Layer Posets}

In this subsection we will prove that a class of posets covering the vast majority of all posets is diametrally reversing.

\begin{figure}[ht]
   \begin{center}
    \psfrag{A}{$A$}
   \psfrag{B}{$B$}
   \psfrag{C}{$C$}
    \includegraphics[scale=0.5]{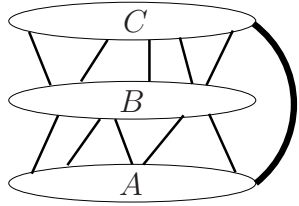}
    \caption{Scheme of a $3$-layer poset. The thick arc indicates that all elements of $A$ are smaller than all elements of $C$.}
    \label{3layer}
   \end{center}
\end{figure}

\begin{definition}
\label{graded}
	A poset~$\P$ is \emph{graded} if every maximal chain of ~$\P$ has the same length. A \emph{$3$-layer poset} is a graded poset of height~$3$ in which each minimum is smaller than each maximum (cf.~Figure~\ref{3layer}).
\end{definition}

\begin{theorem}
\label{3-layer}
Every $3$-layer poset which is not a chain is diametrally reversing. 
\end{theorem}
\proof
Let $\P$ be a 3-layered poset consisting of three layers: the layer of the minima, denoted $A$, the middle layer, denoted $B$, and the layer of the maxima, $C$. First we consider the easy case where there is only one element~$b$ in the middle layer. Now since $\P$ is graded, every minimum has a successor in the middle layer, in this case all minima are smaller than $b$. But this means in fact that every pair of minima is a twin, and by Corollary~\ref{twins} we are done. Hence we can assume $|B|\geq 2$. 

Observe that, if we have an incomparable pair $(a,b)$ with $a \in A$ and $b \in B$, then it forms a critical pair automatically: $\downarrow \!\! a = \emptyset \subseteq \downarrow \!\! b$, and $\uparrow \!\! b \subseteq C \subseteq \uparrow \!\! a$. An analogous argument shows that every incomparable pair $(b,c)$ with $b \in B$ and $c \in C$ is critical. Now let us call an LE \emph{mixing} if it contains an element of $B$ appearing before an element of $A$ or an element of $C$ appearing before an element of $B$. Then we have shown that any mixing LE is automatically reversing. 

So let us consider a non-mixing diametral LE $L$, forming a diametral pair with $L'$. We will first analyse $L'$. If $L'$ is a non-mixing LE, too, then $\P$ must be a complete layered poset, since in all other cases we can find two LEs with larger distance. But every complete layered poset contains many twins, and so by Corollary \ref{twins} we are done. Thus we can assume that $L'$ is a mixing LE. 

The LE $L'$ starts with an initial segment consisting of only minima, followed by a part where elements of the three layers are mixed, and it finishes with a final segment consisting of only maxima. We call the set of elements forming the initial segment $A' \subseteq A$, and the set forming the final segment $C' \subseteq C$. Note that since $\P$ is graded the partition of the elements into three layers is unique, in particular, $A'$ and $C'$ cannot be empty. 

We label the elements of the three layers as $a_1, a_2, \ldots a_k$, $b_1, b_2, \ldots b_{\ell}$ and $c_1, c_2, \ldots c_k$ in the order in which they appear in $L'$. Now we return to $L$. Since $L$ is non-mixing, the elements of every level appear successively. So we are free to choose the layer orders in $L$ and therefore choose them in a way which differs most from $L'$. Hence the elements of each layer appear in $L$ exactly in the opposite order of their order in $L'$, and thus $L$ has the form 
\[
L=a_k a_{k-1} \ldots a_2 a_1 b_{\ell} b_{\ell-1}  \ldots b_2 b_1 c_m c_{m-1} \ldots c_2 c_1.
\]
Assuming $A' = \{a_1,  a_2,  \ldots,  a_{k'}\}$ and $C'=\{c_{m'},  c_{m'+1},  \ldots,  c_m\}$ we know that $L'$ looks like this:
\[
	L'=a_1 a_2 \ldots a_{k'} b_1 X b_{\ell} c_{m'} \ldots c_{m-1} c_m.
\]
Here, $X$ denotes the mixed part of $L'$ consisting of elements of potentially all three layers.

We claim that $(b_1, b_{\ell})$ is a critical pair. Since they are both elements of the middle layer, they are incomparable. The predecessors of $b_1$ are contained in $A'$, since these are the only elements smaller than $b_1$ in $L'$. We want to show that $A' \subseteq \: \downarrow \!\! b_{\ell}$. Now in $L$, the set $A'$ is found at the end of the sequence of minima, immediately before the element $b_{\ell}$. Hence if $a_1$ and $b_{\ell}$ were incomparable, they could be exchanged in $L$ to yield a pair of LEs with larger distance, a contradiction. So we have $a_1 < b_{\ell}$. But in fact we can extend this argument to show that any $a_i \in A'$ needs to be smaller than $b_{\ell}$. This is because we can choose any order of the initial segment $A'$ in $L'$. So we are free to turn any $a_i \in A'$ into $a_1$ by reordering the elements of $A'$ in $L'$ and $L$, and thus with the above argument $a_i < b_{\ell}$ for any $a_i \in A'$. Hence we have $\downarrow \! b_1 \subseteq  A' \subseteq \: \downarrow \! b_{\ell}$. 

It remains to show that $\uparrow \!\! b_{\ell} \subseteq \: \uparrow \!\! b_1$. The argument works analogously: From the position of $b_{\ell}$ in $L'$ we deduce that $\uparrow \! b_{\ell} \subseteq C'$. In~$L$, the element~$b_1$ appears immediately before the elements of $C'$, and if some $c_i \in C'$ were incomparable with $b_1$, then we could construct a pair of LEs with larger distance by rearranging $C'$ in $L$ and $L'$ and exchanging $b_1$ with $c_i$. This would be a contradiction to $L, L'$ being a diametral pair, and thus $\uparrow \!\! b_{\ell} \subseteq \: \uparrow \!\! b_1$. 

We have shown that $(b_1, b_{\ell})$ is a critical pair. In $L$, we have $b_{\ell} < b_1$, so we found a critical pair that is reversed in $L$. This shows that if a non-mixing LE of~$\P$ is diametral, then it is reversing. 
\qed

With the same ideas it can be proved that posets of height 2, which can be thought of as 3-layer posets with $C=\emptyset$, are diametrally reversing. The proof of Theorem~\ref{3-layer} goes through with in fact many steps becoming easier because the empty set is contained in any set.

\begin{proposition}
Every poset of height $2$ which is not a chain is diametrally reversing. 
\end{proposition}

Kleitman and Rothschild~\cite{kleitman-roth-3layer} showed that almost all posets are 3-layer posets, that is, when considering posets on $n$ elements, the proportion of the number of 3-layer posets to the number of all posets tends to 1 as $n$ tends to infinity. Since there is only one chain for each $n$, we can neglect these, and so with Theorem~\ref{3-layer} we immediately obtain the following result: 

\begin{corollary}
Almost all posets are diametrally reversing.
\end{corollary}

\section{Concluding Remarks}
\label{concluding}

There are still open questions about when certain classes of posets are diametrally reversing. For example, are all posets of height 3 diametrally reversing? 

One case of particular interest is that of the Boolean lattices, for which we still know surprisingly little. It seems that one should be able to find out what diametral pairs of LEs of these posets look like. A natural candidate for a diametral pair is formed by the reverse lexicographic LE and the reverse antilexicographic LE. The reverse lexicographic LE $L_1$ takes $1,2,\ldots, n$ as priority order and covers all subsets of $[n]$ formed by an initial segment of it before moving to subsets involving the next element. The reverse antilexicographic LE $L_2$ uses $n, n-1, \ldots,1$ as priority order for this principle. The distance between $L_1$ and $L_2$ is $2^{2n-2}-(n+1)2^{n-1}$. 
Felsner and Reuter~\cite{felsner-led} conjecture that $\textmd{led}(B_n)=2^{2n-2}-(n+1)2^{n-1}$, i.e.~that $L_1$ and $L_2$ form a diametral pair of $B_n$ for all $n$. They prove this for $n \leq 4$. Moreover, they prove that this conjecture would follow if every diametral pair of the Boolean lattice contains a reversing LE. This was in fact their motivation for the conjecture we disproved. 

Are all Boolean lattices diametrally reversing? It seems likely, but small numbers might fool us. If so, what is the smallest $n$ such that $B_n$ is not diametrally reversing? 

As for the complexity question, it seems clear that large width of a poset makes determining its exact LE-diameter more difficult. We think that there probably is a polynomial time algorithm for computing the LE-diameter of posets of arbitrary fixed width, but so far, this question remains open.

\ni {\em Acknowledgement.} We thank Stefan Felsner for fruitful discussions.

\bibliography{my_library}

\begin{thebibliography}{10}

\bibitem{felsner-dilworth}
Stefan Felsner, Vijay Raghavan, and Jeremy Spinrad.
\newblock Recognition algorithms for orders of small width and graphs of small
  dilworth number.
\newblock {\em Order}, 20(4):351--364, 2003.

\bibitem{felsner-led}
Stefan Felsner and Klaus Reuter.
\newblock The linear extension diameter of a poset.
\newblock {\em SIAM J. Disc. Math.}, 12(3):360--373, 1999.

\bibitem{gareyjohnson}
Michael~R. Garey and David~S. Johnson.
\newblock {\em Computers and intractability: A guide to the theory of
  {NP}-completeness}.
\newblock Freeman, 1979.

\bibitem{kleitman-roth-3layer}
Daniel~J. Kleitman and Bruce~L. Rothschild.
\newblock Asymptotic enumeration of partial orders on a finite set.
\newblock {\em Trans. Amer. Math. Soc.}, 205:205--220, 1975.

\bibitem{moehring-tractable}
Rolf~H. M\"ohring.
\newblock Computationally tractable classes of ordered sets.
\newblock In Ivan Rival, editor, {\em Algorithms and Order}, volume 255 of {\em
  NATO ASI Series C}. Kluwer Academic Publishers, 1989.

\bibitem{naatz}
Michael Naatz.
\newblock The graph of linear extensions revisited.
\newblock {\em SIAM J. Disc. Math.}, 13(3):354--369, 2000.

\bibitem{naatz-diss}
Michael Naatz.
\newblock {\em Acyclic Orientations of Mixed Graphs}.
\newblock PhD thesis, Technische Universitat Berlin, 2001.

\bibitem{pruesse91}
Gara Pruesse and Frank Ruskey.
\newblock Generating the linear extensions of certain posets by transpositions.
\newblock {\em SIAM J. Discrete Mathematics}, 4:413--422, 1991.

\bibitem{reuter-comp}
Klaus Reuter.
\newblock The comparability graph and the graph of linear extensions of a
  poset.
\newblock {\em Hamburger Beitr\"age zur Mathematik}, Heft 57, 1996.

\bibitem{reuter}
Klaus Reuter.
\newblock Linear extensions of a poset as abstract convex sets.
\newblock {\em Hamburger Beitr\"age zur Mathematik}, Heft 56, 1996.

\bibitem{ruskey93}
Frank Ruskey and Carla~D. Savage.
\newblock Hamilton cycles that extend transposition matchings in {C}ayley
  graphs of {$S_n$}.
\newblock {\em SIAM J. Discr. Math.}, 6(1):152--166, 1993.

\bibitem{trotter92}
William~T. Trotter.
\newblock {\em Combinatorics and partially ordered sets}.
\newblock The Johns Hopkins University Press, 1992.

\bibitem{west93}
Douglas~B. West.
\newblock Generating linear extensions by adjacent transpositions.
\newblock {\em J. Comb. Theory Ser. B}, 58(1):58--64, 1993.

\end{thebibliography}
\bibliographystyle{plain}
\end{document}